\documentclass[11pt]{amsart}

\usepackage{amsmath}
\usepackage{amssymb}
\usepackage[all]{xy}

\newtheorem{theorem}{Theorem}[section]
\newtheorem{claim}[theorem]{Claim}

\newtheorem{lemma}[theorem]{Lemma}
\newtheorem{facts}[theorem]{Facts}

\newtheorem{proposition}[theorem]{Proposition}
\newtheorem{corollary}[theorem]{Corollary}

\theoremstyle{definition}
\newtheorem{definition}[theorem]{Definition}

\newtheorem{question}[theorem]{Question}

\theoremstyle{remark}
\newtheorem{remark}[theorem]{Remark}

\newcount\skewfactor
\def\mathunderaccent#1#2 {\let\theaccent#1\skewfactor#2
	\mathpalette\putaccentunder}
\def\putaccentunder#1#2{\oalign{$#1#2$\crcr\hidewidth
		\vbox to.2ex{\hbox{$#1\skew\skewfactor\theaccent{}$}\vss}\hidewidth}}
\def\name{\mathunderaccent\tilde-3 }

\def\smallbox#1{\leavevmode\thinspace\hbox{\vrule\vtop{\vbox
			{\hrule\kern1pt\hbox{\vphantom{\tt/}\thinspace{\tt#1}\thinspace}}
			\kern1pt\hrule}\vrule}\thinspace}

\newcommand{\cf}{{\rm cf}}

\newcommand{\dom}{{\rm dom}}
\newcommand{\stem}{{\rm stem}}
\newcommand{\len}{{\rm len}}

\newcommand{\la}{\langle}
\newcommand{\ra}{\rangle}
\newcommand{\fr}{{}^\frown}
\newcommand{\uhr}{\upharpoonright}
\newcommand{\force}{\Vdash}
\newcommand{\can}{\check}
\newcommand{\power}{\mathcal{P}}
\newcommand{\MS}{\mathcal{MS}}

\def\qedref#1{$\qed_{\reforiginal{#1}}$}

\title{Weak prediction principles}

\author{Omer Ben-Neria}
\address{Department of Mathematics, University of California, Los Angeles}
\email{obneria@math.ucla.edu}

\author{Shimon Garti}
\address{Institute of Mathematics,
	The Hebrew University of Jerusalem,
	Jerusalem 91904, Israel}
\email{shimon.garty@mail.huji.ac.il}

\author{Yair Hayut}
\address{Institute of Mathematics,
	The Hebrew University of Jerusalem,
	Jerusalem 91904, Israel}
\email{yair.hayut@mail.huji.ac.il}

\subjclass[2010]{03E05}
\keywords{Weak diamond, very weak diamond, Radin forcing}

\begin{document}
	\let\labeloriginal\label
	\let\reforiginal\ref

	\begin{abstract}
		We prove the consistency of the failure of the weak diamond $\Phi_\lambda$ at strongly inaccessible cardinals. On the other hand, we show that the very weak diamond $\Psi_\lambda$ is equivalent to the statement $2^{<\lambda}<2^\lambda$ and hence holds at every strongly inaccessible cardinal.
	\end{abstract}
	
	\maketitle
	
	\newpage
	
	\section{Introduction}
	
	The prediction principle $\Diamond_\lambda$ (diamond on $\lambda$) was discovered by Jensen, \cite{MR0309729}, who proved that it holds over any regular cardinal $\lambda$ in the constructible universe. This principle says that there exists a sequence $\langle A_\alpha:\alpha<\lambda\rangle$ of sets, $A_\alpha\subseteq\alpha$ for every $\alpha<\lambda$, such that for every $A\subseteq\lambda$ the set $\{\alpha<\lambda:A\cap\alpha=A_\alpha\}$ is stationary.
	
	Jensen introduced the diamond in 1972, and the main focus was the case of $\lambda=\aleph_1$. It is immediate that $\Diamond_{\aleph_1}\Rightarrow 2^{\aleph_0}=\aleph_1$, but consistent that $2^{\aleph_0}=\aleph_1$ along with $\neg\Diamond_{\aleph_1}$. Motivated by algebraic constructions, Devlin and Shelah \cite{MR0469756} introduced a weak form of the diamond principle which follows from the continuum hypothesis:
	
	\begin{definition}[The Devlin-Shelah weak diamond]
		\label{wweakdiamond}
		Let $\lambda$ be a regular uncountable cardinal. \newline
		The weak diamond on $\lambda$ (denoted by $\Phi_\lambda$) is the following principle: \newline
		For every function $c:{}^{<\lambda} 2\rightarrow 2$ there exists a function $g\in{}^{\lambda} 2$ such that $\{\alpha\in\lambda: c(f\upharpoonright \alpha)=g(\alpha)\}$ is a stationary subset of $\lambda$ whenever $f\in{}^{\lambda} 2$.
	\end{definition}
	
	The idea is that we replace the prediction of the initial segments of a set (or a function) by predicting only their color. The function $c$ is a coloring, and the function $g$ is the weak diamond function which gives stationarily many guesses for the $c$-color of the initial segments of every function $f$. It is easy to see that the real diamond implies the weak diamond.
	
	Concerning cardinal arithmetic, $\Phi_{\aleph_1}$ follows indeed from the continuum hypothesis. Moreover, the weak continuum hypothesis $2^{\aleph_0}<2^{\aleph_1}$ implies $\Phi_{\aleph_1}$ as proved in \cite{MR0469756}. On the other hand, $\Phi_{\aleph_1}$ implies $2^{\aleph_0}<2^{\aleph_1}$ (as noted by Uri Abraham) so the two assertions are equivalent.
	
	The diamond and the weak diamond are prediction principles, while $2^{\aleph_0}<2^{\aleph_1}$ and $2^{\aleph_0}=\aleph_1$ belong to cardinal arithmetic. The above statements show that a simple connection may exist. Actually, if $\lambda>\aleph_1$ and $\lambda=\kappa^+$ then the situation becomes even simpler:
	
	\begin{proposition}[Diamonds and cardinal arithmetic]
		\label{prpreferee} Assume that $\lambda=\kappa^+>\aleph_1$.
		Then $\Diamond_\lambda\Leftrightarrow 2^\kappa=\kappa^+$ and $\Phi_\lambda \Leftrightarrow 2^\kappa<2^{\kappa^+}$.
	\end{proposition}

	For the assertion concerning the diamond principle see \cite{MR2596054}. The assertion for the weak diamond is a straightforward generalization of \cite{MR0469756}, the easy direction appears explicitly in \cite{MR3604115} Proposition 1.2 and the substantial direction can be extracted from \cite{MR0469756} upon replacing $\aleph_1$ by $\kappa^+$. The negation of $\Diamond_{\aleph_1}$ with the affirmation of $2^{\aleph_0}=\aleph_1$ has been recognized as a peculiarity of the first uncountable cardinal.
	
	However, the situation is totally different if $\lambda=\cf(\lambda)$ is a limit cardinal. One direction is still easy and has a general nature. If $\Diamond_\lambda$ then $2^{<\lambda}=\lambda$ and if $\Phi_\lambda$ then $2^{<\lambda}<2^\lambda$ (see Claim \ref{mc} below). Notice that these formulations coincide with the above description of the successor case $\lambda=\kappa^+$. But can we prove the opposite implication?
	
	If $\lambda=\cf(\lambda)$ is a limit cardinal then $\lambda$ is a large cardinal
	(in the philosophical sense; its existence cannot be established in ZFC).
	Now if $\lambda$ is a large enough cardinal, such as measurable, then $\Diamond_\lambda$ holds.
	Furthermore, if $\lambda$ is ineffable (a property which may live happily with ${\rm V}={\rm L}$)
	or even subtle then $\Diamond_\lambda$ and hence also $\Phi_\lambda$ \cite{JensenKunan}.
	
	A striking (and unpublished) result of Woodin shows that relatively small large cardinals
	are different than successor cardinals. Woodin proved the consistency of a strongly inaccessible cardinal $\lambda$ for which $\Diamond_\lambda$ fails. Moreover, the construction can be strengthened to strongly Mahlo.
	In the current paper we prove the same assertion, upon replacing the diamond by the weak diamond. For both the diamond and the weak diamond we do not know what happens if $\lambda$ is weakly compact.
	
	We conclude that the cardinal arithmetic assumption $2^{<\lambda}<2^\lambda$ is strictly weaker than the prediction principle $\Phi_\lambda$. It is still tempting to look for a prediction principle which is characterized by $2^{<\lambda}<2^\lambda$. We define the following:
	
	\begin{definition}[The very weak diamond]
		\label{vwddef}
		Let $\lambda$ be an uncountable cardinal. \newline
		The very weak diamond on $\lambda$ (denoted by $\Psi_\lambda$) is the following principle: \newline
		For every function $c:{}^{<\lambda} 2\rightarrow 2$ there exists a function $g\in{}^{\lambda} 2$ such that $\{\alpha\in\lambda: c(f\upharpoonright \alpha)=g(\alpha)\}$ is an unbounded subset of $\lambda$ whenever $f\in{}^{\lambda} 2$.
	\end{definition}
	
	As we shall see, $\Psi_\lambda\Leftrightarrow 2^{<\lambda}<2^\lambda$ whenever
	$\lambda=\cf(\lambda)>\aleph_0$. Let us mention Galvin's property which
	says that every collection $\{C_\alpha:\alpha<\lambda\}$ of club subsets of $\kappa^+$ has a sub-collection of size $\kappa^+$ whose intersection is a club.
	This principle follows
	from $2^{<\lambda}<2^\lambda$ but consistent with $2^{<\lambda}=2^\lambda$ (see \cite{MR3604115}). Bringing all these principles together we have a systematic hierarchy of weak prediction principles, each of which implied by the stronger one but strictly weaker from the next stage.
	
	The basic tool for proving $\neg\Phi_\lambda$ over small large cardinals is Radin forcing, \cite{MR670992}. Since there are several ways to introduce this forcing notion we indicate that our approach is taken from \cite{MR2768695}, and in particular we use the Jerusalem notation, i.e., $p\leq_{\mathbb{R}}q$ means that $q$ is stronger than $p$.
	
	\newpage
	
	\section{Very weak diamond and cardinal arithmetics}
	We commence with the easy direction about the connection between $\Psi_\lambda$ and cardinal arithmetic.
	
	We commence with the easy direction about the connection between $\Phi_\lambda$ and cardinal arithmetic. Actually, the claim below applies to $\Psi_\lambda$ as well. A parallel assertion can be proved easily for $\Diamond_\lambda$.
	
	\begin{theorem}[The basic claim]
		\label{mc}
		Let $\lambda$ be an uncountable cardinal. \newline
		If $\Psi_\lambda$ holds then $2^{<\lambda}<2^\lambda$.
	\end{theorem}
	
	\par\noindent\emph{Proof}. \newline
	Assume that $2^{<\lambda}=2^\lambda$. We will show that $\Psi_\lambda$ fails. Let $b:2^{<\lambda}\rightarrow 2^\lambda$ be a surjection, such that for every $\delta < \alpha < \lambda$ and every $t\in{}^\alpha 2$, if $\varepsilon\in[\delta,\alpha)\Rightarrow t(\varepsilon)=0$, then $b(t)=b(t\upharpoonright\delta)$.
	We are trying to describe a coloring $F:{}^{<\lambda}2\rightarrow 2$, which exemplifies the failure of the very weak diamond.
	
	Assume $\alpha<\lambda$ and $\eta\in{}^\alpha 2$. Set $F(\eta)=[b(\eta)](\alpha)$.
	Let $g\in{}^\lambda 2$ be a function. We will show that $g$ does not predict $F$. Let $h\in{}^\lambda 2$ be the opposite function, i.e. $h(\alpha)=1-g(\alpha)$ for every $\alpha<\lambda$. Let $t\in{}^{<\lambda}2$ be any mapping for which $b(t)=h$. Let $\delta<\lambda$ be such that $t\in{}^\delta 2$. We define $f\in{}^\lambda 2$ as an extension of $t$ as follows. If $\alpha<\delta$ then $f(\alpha)=t(\alpha)$ and if $\alpha\geq\delta$ then $f(\alpha)=0$. Observe that $b(f\upharpoonright\alpha)=b(f\upharpoonright\delta)$ for every $\alpha\in(\delta,\lambda)$. Consequently, \[F(f\upharpoonright\alpha)=[b(f\upharpoonright\alpha)](\alpha) =[b(f\upharpoonright\delta)](\alpha)=[b(t)](\alpha) =h(\alpha)\neq g(\alpha),\] so $g$ fails to predict $F(f\upharpoonright\alpha)$ on an end-segment of $\lambda$, as wanted.
	\hfill \qedref{mc}
	It follows from the above claim that the weak diamond and the very weak diamond are equivalent in the successor case. We comment that the same holds for the common diamond over successor cardinals, if one replaces the requirement of stationary set of guesses by an unbounded set.
	That is, if $\Diamond'_\lambda$ says that every $A\subseteq\lambda$ is guessed by an unbounded set of $A_\alpha$'s then $\Diamond_\lambda\Leftrightarrow\Diamond'_\lambda$ whenever $\lambda=\kappa^+$ (see \cite{Devlin} for the case of $\aleph_1$. The same argument works for the general successor case). For the following corollary we recall that a cardinal $\lambda$ is weakly inaccessible if and only if $\lambda$ is a regular limit cardinal.
	
	\begin{corollary}[$\Phi_\lambda$ and $\Psi_\lambda$]
		\label{cor1}
		If $\lambda=\kappa^+$ then $\Phi_\lambda\Leftrightarrow\Psi_\lambda$. The same holds true if $\lambda$ is weakly inaccessible and $2^{<\lambda}=2^\kappa$ for some $\kappa<\lambda$.
	\end{corollary}
	
	\par\noindent\emph{Proof}. \newline
	The implication $\Phi_\lambda\Rightarrow\Psi_\lambda$ results from the definition
	for both cases mentioned in the above statement.
	For the opposite direction, if $\Psi_\lambda$ then $2^{<\lambda}<2^\lambda$ by the above theorem, and hence $\Phi_\lambda$ follows from Proposition \ref{prpreferee} in the successor case and from Claim \ref{mmc} in the weakly inaccessible case.
	
	\hfill \qedref{cor1}
	
	Back to the general case, Theorem \ref{mc} gives one direction by showing that $\Phi_\lambda$ implies $2^{<\lambda}<2^\lambda$.
	As we shall see in Theorem \ref{mt} below, the other direction cannot be proved. Indeed, every strongly inaccessible cardinal satisfies $2^{<\lambda}<2^\lambda$, but the negation of $\Phi_\lambda$ can be forced over some strongly inaccessible cardinal. However, we can prove that $2^{<\lambda}<2^\lambda$ is equivalent to the very weak diamond $\Psi_\lambda$. We shall use the fact that if $\lambda$ is weakly inaccessible and $2^{<\lambda}=2^\kappa<2^\lambda$ for some $\kappa<\lambda$ then $\Phi_\lambda$. This fact appears already in \cite{MR0469756}, without explicit proof. Since it plays a key-role in the theorem below, we spell out the proof:
	
	\begin{claim}[Weak diamond out of nowhere]
		\label{mmc}
		If $\lambda$ is weakly inaccessible, $2^{<\lambda}<2^\lambda$ and $2^\kappa=2^{<\lambda}$ for some $\kappa<\lambda$, then $\Phi_\lambda$.
	\end{claim}
	
	\par\noindent\emph{Proof}. \newline
	Assume towards contradiction that $\neg(\Phi_\lambda)$, and choose a coloring $F:{}^{<\lambda}2\rightarrow 2$ which exemplifies it.
	Denote the club filter over $\lambda$ by $\mathcal{D}_\lambda$.
	
	Let $\mathcal{S}$ be the collection of all the sequences $s$, of length $1 + \beta  + \beta$ such that:
	\begin{itemize}
		\item $s(0), \beta < \lambda$.
		\item $s(1 + i), s(1 + \beta + i) \in {}^{s(0)} 2$ for all $i < \beta$.
	\end{itemize}
	We will denote $\alpha^s = s(0)$, $g^s_\nu = s(1 + \nu)$, $f^s_\nu = s(1 + \beta + \nu)$. We will omit the superscript $s$ where it is clear from the context.
	
	Observe that $|\mathcal{S}|=2^{<\lambda}$, which equals $2^\kappa$ by the assumption of the theorem.
	Let $h$ be a one-to-one mapping between $\mathcal{S}$ and ${}^\kappa 2$.
	
	Assume $d\in{}^\lambda 2$ is any function. We will define, by induction on $n < \omega$, functions $g_{\kappa \cdot n + \delta}$, $f_{\kappa \cdot n + \delta}$, $\delta < \kappa$, and clubs $C_n$ so that each $f_\eta$ codes $g_\eta$ using the failure of $\Phi_\lambda$. We will show, eventually, that this coding process produces an unique code from each function from ${}^\lambda 2$ and thus enables us to obtain a surjection from ${}^{{<}\lambda}2$ onto ${}^\lambda 2$.
	
	For $n=0$ let $g_0 = d, f_0=f$ for some $f$ such that the set $\{\alpha \mid F(f\restriction \alpha) = g(\alpha)\}$ contains a club, $C$. We set $C_0=C$. Let $g_\delta, f_\delta$ be arbitrary for $0 < \delta < \kappa$.
	
	Let $n > 0$. Let us assume that $f_\nu, g_\nu$ where defined for every $\nu < \kappa \cdot n$ and let us assume that $C_{n - 1}$ is defined. We want to define the functions $f_{\kappa \cdot n + \delta}, g_{\kappa \cdot n + \delta}$ and the club $C_n$.
	
	For each $\alpha < \lambda$ let $\beta_{\alpha,\eta}$ be the first member of $C_{n - 1}$ greater than $\alpha$. We define simultaneously the sequence of functions $\langle g_{\kappa\cdot n + \delta} \mid \delta<\kappa\rangle$. For this end, we have to determine the value of all these functions for every $\alpha<\lambda$.
	
	For $\alpha < \lambda$, let $s^n_\alpha$ be the following member of $\mathcal{S}$:
	\begin{itemize}
		\item $s^n_\alpha(0) = \beta_{\alpha,\eta}$,
		\item $s^n_\alpha(1 + \nu) = g_\nu\upharpoonright\beta_{\alpha,\eta}$ for every $\nu < \kappa \cdot n$ and
		\item $s^n_\alpha(1 + \beta_{\alpha,\eta} + \nu) = f_\nu\upharpoonright\beta_{\alpha,\eta}$ for every $\nu < \kappa \cdot n$. \end{itemize}
	Recall that $h(s^n_\alpha) \in {}^\kappa 2$. Let us define $g_{\kappa \cdot n + \delta}(\alpha) = h(s^n_\alpha)(\delta)$.
	
	Having the functions $g_{\kappa\cdot n +\delta}$ at hand for every $\delta<\kappa$, we choose for each one of them a function $f_{\kappa\cdot n +\delta}\in{}^\lambda 2$ such that \[A_{\kappa\cdot n +\delta} = \{\alpha<\lambda: F(f_{\kappa\cdot n +\delta} \upharpoonright\alpha)=g_{\kappa\cdot n +\delta}(\alpha)\}\in \mathcal{D}_\lambda.\] Finally, we choose a club $C_{n}$ of $\lambda$ such that \[C_{n}\subseteq \bigcap\limits_{\delta<\kappa}A_{\kappa\cdot n +\delta}\cap C_{n - 1}.\]
	
	Let $\gamma = \min \bigcap_{n < \omega} C_n$. Let us define the \emph{code} of $d$, $\mathcal{C}^d$, to be the following sequence of length $\kappa \cdot \omega \cdot 2$:
	\begin{enumerate}
		\item $\mathcal{C}^d(0) = \gamma$.
		\item $\mathcal{C}^d(1 + \nu) = g_\nu \restriction \gamma$ for $\nu < \kappa \cdot \omega$.
		\item $\mathcal{C}^d(\kappa\cdot \omega + \nu) = f_\nu\restriction \gamma$ for $\nu < \kappa  \cdot \omega$.
	\end{enumerate}
	
	$\mathcal{C}^d \in \lambda \times {}^{\kappa \cdot \omega \cdot 2} ({}^{{<}\lambda}2)$, so there are $2^{{<}\lambda}$ possible values for the code $\mathcal{C}^d$. Let us show that one can reconstruct $d$ from $\mathcal{C}^d$ and thus conclude that $2^{{<}\lambda} = 2^{\lambda}$.
	
	Indeed, let $\mathcal{C}^d$ be the code of some function $d$. Let $g_\nu, f_\nu$ and $C_n$ be the functions and the clubs which were generated in the course of the construction of $\mathcal{C}^d$. Let $d'$ be another function and let $g'_\nu, f'_\nu, C'_n$ be the functions and the clubs which were generated in the course of constructing $\mathcal{C}^{d'}$. Let us assume that $\mathcal{C}^d = \mathcal{C}^{d'}$ and show that $d = d'$.
	
	Let us show, by induction on $\alpha \in \bigcap C_n$, that $g_\nu \restriction \alpha = g'_\nu \restriction \alpha$, $f_\nu\restriction \alpha = f'_\nu \restriction \alpha$ and $\bigcap_{n < \omega} C_n \cap (\alpha + 1) = \bigcap_{n < \omega} C'_n \cap (\alpha + 1)$. This is enough, as $g_0 = d$, $g_0' = d'$.
	
	For simplicity of notations, let $D = \bigcap_{n < \omega} C_n$, $D' = \bigcap_{n < \omega} C'_n$.
	
	For $\alpha = \min D = \min D'$, the inductive assumption holds since $\mathcal{C}^d = \mathcal{C}^{d'}$.
	
	Let us assume that the claim is true for every $\beta < \alpha$ in $D$ and let us show its validity for $\alpha$. If $\alpha$ is an accumulation point of $D$ then $\alpha \in D$, since $D$ is a club. Similarly, $\alpha \in D'$. The rest of the inductive assumption holds trivially.
	
	Let $\alpha$ be non-accumulation point, above the minimal point of $D$, and let $\gamma = \max D \cap \alpha$. Since the clubs $C_n$ are decreasing, $\sup \beta_{n, \gamma} = \alpha$.
	
	Since $\gamma \in C_n \cap C'_n$ for all $n$ and since $f_\nu \restriction \gamma = f'_\mu \restriction \gamma$ for all $\nu < \kappa \cdot \omega$,
	\[F(f_{\nu} \restriction \gamma) = g_{\nu}(\gamma) = F(f'_{\nu} \restriction \gamma) = g'_{\nu}(\gamma)\]
	But $g_{\kappa \cdot n + \delta}(\gamma) = h(s^n_\gamma)(\delta)$ and since $h$ is one to one, we conclude that $\beta_{n,\gamma} = \beta'_{n, \gamma}$, $g_\nu \restriction \beta_{n, \gamma} = g'_\nu \restriction \beta_{n, \gamma}$ and $f_\nu \restriction \beta_{n, \gamma} = f'_\nu \restriction \beta_{n, \gamma}$, for all $\mu < \kappa \cdot n$. This is true for all $n < \omega$, and thus we conclude that the induction assumption holds for $\alpha$.
	\hfill \qedref{mmc}
	
	\begin{remark}
		\label{r} Another way to phrase the idea in the above proof is by noticing that if $2^{<\lambda}=2^\kappa$ then $\neg\Phi_\lambda$ codes a one-to-one mapping from $2^\lambda$ into $2^{<\lambda}$, and hence $2^{<\lambda}=2^\lambda$.
	\end{remark}
	
	\hfill \qedref{r}
	
	Now we can prove the following:
	
	\begin{theorem}[Very weak diamond and cardinal arithmetic]
		\label{mmt}
		For every regular uncountable cardinal $\lambda$ we have $2^{<\lambda}<2^\lambda$ iff $\Psi_\lambda$.
	\end{theorem}
	
	\par\noindent\emph{Proof}. \newline
	The successor cardinal case follows from Theorem \ref{mc} and the results of \cite{MR0469756}, so we may assume that $\lambda$ is a limit cardinal.
	If $\Psi_\lambda$ holds then $2^{<\lambda}<2^\lambda$ holds by Theorem \ref{mc}. For the opposite direction, assume that $2^{<\lambda}<2^\lambda$. If $2^{<\lambda}=2^\kappa=2^\lambda$ for some $\kappa<\lambda$ then $\Phi_\lambda$ holds by Claim \ref{mmc} and hence also $\Psi_\lambda$. So assume that this is not the situation (as always happens in the case of a strongly inaccessible cardinal).
	
	We claim that $\Phi_\alpha$ holds for unbounded set of $\alpha$'s below $\lambda$. For proving this assertion, define $C = \{\kappa<\lambda: \forall\gamma<\kappa, 2^\gamma<2^\kappa\}$. By the assumptions on $\lambda$, $C$ is a club subset of $\lambda$. Enumerate the uncountable elements of $C$ by $\{\alpha_\varepsilon:\varepsilon<\lambda\}$.
	
	We claim that $\Phi_{\alpha_\varepsilon}$ holds whenever $\varepsilon$ is a successor ordinal. So let $\varepsilon = \zeta+1$. If $\alpha_\varepsilon$ is a successor cardinal then $\Phi_{\alpha_\varepsilon}$ follows from \cite{MR0469756}. If not, then $2^\gamma = 2^\zeta$ for every sufficiently large $\gamma\in[\alpha_\zeta,\alpha_\varepsilon)$. In which case, $\alpha_\varepsilon$ is regular by the Bukovsky-Hechler Theorem (see Corollary 5.17 in \cite{MR1940513} and the historical notes in p. 61 there) and $\Phi_{\alpha_\varepsilon}$ follows from Claim \ref{mmc} with $\alpha_\zeta$ here standing for $\kappa$ there.
	
	We claim now that $\Psi_\lambda$ holds.
	For this, let $c:{}^{<\lambda}2 \rightarrow 2$ be a coloring. For every $\zeta<\lambda$ let $\varepsilon=\zeta+1$ and let $c_\varepsilon$ be the restriction $c\upharpoonright 2^{<\alpha_\varepsilon}$.
	Choose a function $g_\varepsilon$ which exemplifies $\Phi_{\alpha_\varepsilon}$ with respect to $c_\varepsilon$ for every $\varepsilon=\zeta+1<\lambda$. Define $h:\lambda\rightarrow\lambda$ as follows. For every $\beta<\lambda$ let $h(\beta)$ be the first ordinal $\varepsilon<\lambda$ so that $\alpha_\varepsilon\leq\beta<\alpha_{\varepsilon+1}$.
	
	We define $g\in{}^\lambda 2$ as follows.
	Given $\beta<\lambda$, if $h(\beta)=\varepsilon$ is a limit ordinal then $g(\beta)=0$. If $h(\beta)=\varepsilon$ is a successor ordinal then $g(\beta)=g_\varepsilon(\beta)$.
	Let us show that $g$ exemplifies $\Psi_\lambda$.
	
	Assume $f\in{}^\lambda 2$. For every successor ordinal $\varepsilon=\zeta+1<\lambda$ let $f_\varepsilon=f\upharpoonright\alpha_\varepsilon$. By $\Phi_{\alpha_\varepsilon}$ we can choose an ordinal $\beta_\varepsilon\in [\alpha_\zeta,\alpha_\varepsilon)$ for which $g_\varepsilon(\beta_\varepsilon)=c_\varepsilon(f_\varepsilon\upharpoonright \beta_\varepsilon)$. By the above definitions it follows that $g(\beta_\varepsilon)=c(f\upharpoonright\beta_\varepsilon)$. Since we have unboundedly many $\beta_\varepsilon$ of this form, we are done.
	\hfill \qedref{mmt}
	\section{Failure of weak diamond at strongly inaccessible cardinal}
	In this section we demonstrate the fact that $\Psi_\lambda$ is strictly weaker than $\Phi_\lambda$.
	
	Our next goal is to demonstrate the fact that $\Psi_\lambda$ is strictly weaker than $\Phi_\lambda$. We shall use Radin forcing $R(\vec{U})$.
	For the most part, our presentation and arguments follow  Gitik's handbook chapter, \cite{MR2768695}.
	
	\begin{definition}[Measure sequences]
		Let $\kappa$ be a cardinal and $\vec{U} = \la \kappa\ra \fr \la U_\alpha \mid \alpha < \ell(\vec{U})\ra$ be a sequence such that each $U_\alpha$ is a measure on $V_\kappa$ (i.e., a $\kappa$ complete normal ultrafilter on $V_\kappa$).
		For each $\beta  < \ell(\vec{U})$, let $\vec{U}\uhr \beta$ denote the initial segment $\la \kappa \ra \fr \la U_\alpha \mid \alpha < \beta\ra$. In particular, $\vec{U} \uhr 0 = \langle \kappa\rangle$.
		
		We say $\vec{U}$ is a \textbf{measure sequence} on $\kappa$ if either $\vec U = \langle \kappa \rangle$ or $\ell(\vec{U}) > 0$ and there exists an elementary embedding $j : V \to M$, whose critical point is $\kappa$ and ${}^\kappa M \subset M$, such that for each $\beta < \ell(\vec{U})$, $\vec{U}\uhr \beta \in M$ and $U_\beta = \{ X \subset V_\kappa \mid \vec{U}\uhr \beta \in j(X)\}$.
	\end{definition}
	
	Note that the length of a measure sequence $\vec{U}$, as sequence of sets, is $1 + \ell(U)$. If $\ell(U) > 0$, let $\cap \vec{U}$ denote the filter $\bigcap_{\alpha  < \ell(\vec{U})}U_\alpha$, and $\MS(\kappa)$ denote the set of measure sequences $\vec{\mu}$ on measurable cardinals below $\kappa$. If $\kappa$ is clear from the context, we omit the subscript.
	
	Let $\vec{\mu}$ be a sequence of the form $\la \nu \ra  \fr \la u_i \mid i < \ell(\vec{\mu})\ra$, where each $u_i$ is a measure on $V_{\nu}$. We denote $\nu$ by $\kappa(\vec{\mu})$, and $\cap \{u_i \mid i < \ell(\vec{\mu})\}$ by $\cap \vec{\mu}$. If $\ell(\vec\mu) = 0$, $\cap \vec\mu$ is undefined.
	
	The following lemma is well known (see \cite{MR2768695}):
	\begin{lemma}
		Let $\vec{U}$ be a measure sequence on $\kappa$ with nonzero length. Let $\langle A_{\vec{\nu}} \mid \nu \in \MS(\kappa)\rangle$ be a sequence of sets such that $A_{\vec\nu} \in \cap \vec{U}$. Then the \textbf{diagonal intersection} defined by:
		
		\[A^\star = \triangle_{\vec\nu} A_{\vec\nu} = \{\vec\mu \in \MS(\kappa) \mid \forall \vec\nu \in \MS(\kappa) \cap V_{\kappa(\vec\mu)},\,\vec\mu \in A_{\vec\nu}\}\]
		belongs to $\cap \vec{U}$.
	\end{lemma}
	
	Clearly, if one takes a diagonal intersection over smaller set than $\MS(\kappa)$, the resulting set is only larger and thus belongs to $\cap \vec{U}$.
	
	Let us fix a measurable cardinal $\kappa$ and let $\vec{U}$ be a measure sequence on $\kappa$ of nonzero length. Before we can define the Radin forcing $R(\vec{U})$, we need the to introduce the following sequence of sets $A^n \subset \MS(\kappa)$, $n < \omega$. Let
	$A^0 = \MS(\kappa)$, and  for each $n < \omega$, set $A^{n+1} = \{\vec{\mu} \in A^{n} \mid A^n \cap V_{\kappa(\vec{\mu})} \in \cap \vec{\mu}\}$. We finally define \[\bar{A} = \bigcap_{n<\omega} A^n.\] Since each $A^n$ belongs to $\cap \vec{U}$,  $\bar{A} \in \cap \vec{U}$ as well.
	
	\begin{definition}[Radin forcing]
		The Radin partial ordered set, $R(\vec{U})$, consists of all finite sequences $p = \la d_i \mid i \leq k\ra$ satisfying the following conditions.
		\begin{enumerate}
			\item [$(\aleph)$] $\vec{d} = \la d_i \mid i < k\ra$ is a finite sequence. For every $i \leq k$, $d_i$ is either of the form $\la \kappa_i \ra$ where $\kappa_i < \kappa$ is an ordinal,
			or of the form $d_i = \la \vec{\mu_i} , a_i\ra$ where $\vec{\mu_i}$ is a measure sequence on a measurable cardinal $\kappa_i = \kappa(\vec{\mu_i}) \leq \kappa$ and $a_i \in \cap \vec{\mu_i}$. For each $i \leq k$ we denote $\kappa_i$ by $\kappa(d_i)$ and $a_i$ by $a(d_i)$.
			\item [$(\beth)$] $\la \kappa_i \mid i < k\ra$ is increasing.
			\item [$(\gimel)$] $d_k = \la \vec{U}, A\ra$ for some $A \in \bigcap \vec{U}$ which is a subset of $\bar{A}$.
		\end{enumerate}
		
		Given a condition $p = \la d_i \mid i \leq k\ra$ as above, we will frequently separate its top part $\la \vec{U} , A\ra$ from the other components, and write $p = \vec{d}  \fr \la \vec{U}, A\ra$, 
		where $\vec{d} = \la d_i \mid i < k\ra$. 
		We refer to  $\vec{d}$ as the \textbf{stem} of $p$ and denote it by $\stem(p)$. We say that
		a condition $p^* =  \la d^*_i \mid i \leq k^*\ra$ is a \textbf{direct extension} of $p = \la d_i \mid i \leq k\ra$ if $k^* = k$ and $a(d^*_i) \subset a(d_i)$ for all $i \leq k$ for which $a(d_i)$ exists.
		A condition $p'$ is a \textbf{one-point extension} of $p$ if there exists $j \leq k$ and a measure sequence $\vec{\nu} \in a(d_j)$ with $\kappa(\vec{\nu}) > \kappa(d_{j-1})$,
		and $p'$ is either $\la d_i \mid i < j\ra \fr \la \vec{\nu} \ra \fr \la d_i \mid i \geq j\ra$ if $\vec{\nu} = \langle \alpha\rangle$ for some ordinal $\alpha$, or $\la d_i \mid i < j\ra \fr \la \vec{\nu}, a(d_j) \cap V_{\kappa(\vec{\nu})}\ra \fr \la d_i \mid i \geq j\ra$, where $\vec{\nu}$ is a nontrivial measure sequence. We refer to $p'$ as the one-point extension of $p$ by $\vec{\nu}$ and further denote it by $p \fr \la \vec{\nu} \ra$.
		
		We note that our presentation, which follows \cite{MR2768695}, abuses the symbol "$\fr$" in the context of Radin forcing. Therefore, when we use the concatenation symbol, $\fr$, with a member of the Radin forcing as the left argument (i.e., $p \fr \vec{\nu}$) its meaning is the weakest condition in the Radin forcing which is stronger than the left argument (i.e., $p$) and contains the right argument in its sequence (i.e., contains $\vec{\nu}$). This interpretation differs from the standard one which is appending the right argument to the left sequence. 
		
		A condition $q$ extends $p$ if it is obtained from $p$ by a finite sequence of one-point extensions and direct extensions.
	\end{definition}
	
	Let $\vec{U} =  \la \kappa \ra \fr \la U_\tau \mid \tau < \kappa^+\ra$ be a measure sequence of length $\kappa^+$, derived from an elementary embedding $j : V \to M$ as above. The following results are established in  \cite{MR2768695}.
	\begin{facts}
		\begin{enumerate}
			\item $R(\vec{U})$ satisfies $\kappa^+.c.c$.
			\item $R(\vec{U})$ is a Prikry type forcing notion. Namely, for every condition $p \in R(\vec{U})$ and a statement $\sigma$ of the forcing language, $p$ has a direct extension $p^*$ which decides $\sigma$.
			\item $R(\vec{U})$ preserves all cardinals.
			\item $\kappa$ remains regular and strong limit  in a $R(\vec{U})$ generic extension.
			\item Suppose that $G \subset R(\vec{U})$ is a generic filter.
			Let $MS_G \subset \MS$ be the set of all $\vec{u} \in \MS$ for which there exists some $p \in G$ of the form $p = \vec{d} \fr \la \vec{U}, A\ra$\footnote{Here, $\vec{d} \fr \la \vec{U},A\ra$ does not denote a one-point extension but rather the condition which is obtained by appending the top component $\vec{U},A\ra$ to the finite sequence $\vec{d} \in V_\kappa$.} such that $\vec{d} = \la u_0,a_0\ra ,\dots , \la u_{k-1},a_{k-1}\ra$ and $\vec{u} = u_i$ for some $i < k$.
			The set $C_G = \{ \kappa(\vec{u}) \mid \vec{u} \in MS_G\}$ is a closed unbounded subset of $\kappa$, called the generic Radin club associated with $G$.
			\item $V[MS_G] = V[G]$, moreover $G$ is the set of all $p \in R(\vec{U})$ such that $p = \vec{d} \fr \la \vec{U}, A\ra$, $\vec{d} = \la u_0,a_0\ra ,\dots , \la u_{k-1},a_{k-1}\ra$ and $u_i\in MS_G$ for all $i < \len(d)$.
			\item Let $\vec{\nu} \in MS_G$ be a measure sequence of positive length. Let $G(\vec{\nu})$ be the filter in $R(\vec{\nu})$ generated by $MS_G \cap V_\nu$. Then $G(\vec{\nu})$ is a generic filter for $R(\vec{\nu})$.
		\end{enumerate}
	\end{facts}

	\begin{theorem}[Strong inaccessibility and $\neg\Phi_\kappa$]
		\label{mt}
		Assuming the existence of a measure sequence $\vec{U} =\la \kappa\ra \fr \la U_\alpha \mid \alpha < \kappa^+\ra$ derived from an embedding $j  : V \to M$, such that $\vec{U} \in M$ and  $M\models 2^\kappa=2^{\kappa^+}$, it is consistent that there is an inaccessible cardinal $\kappa$ such that $\neg\Phi_\kappa$.
	\end{theorem}
	
	The large cardinal assumption of the Theorem is known to be consistent. For example, if $j : V \to M$ is a $\kappa^+$-supercompact embedding then $M$ is closed under sequences of length $\kappa^+$ and thus contains $\vec{U}$. By a well-known argument of Silver, the existence of such a supercompact embedding $j : V \to M$ is consistent with $2^\kappa = 2^{\kappa^+} = \kappa^{++}$ (see \cite[Section 12]{Cummings-HBchapter}). The supercompact assumption can be further reduced to a hyper measurability assumption which involves a measurable cardinal of some high Mitchell order, using an argument of Woodin (see \cite{Gitik-negSCH}).
	
	\par\noindent\emph{Proof}. \newline
	Suppose $\vec{U} = \la \kappa \ra \fr \la U_\tau \mid \tau < \kappa^+\ra$ is a measure sequence of a measurable cardinal $\kappa$ derived from an elementary embedding $j : V \to M$
	satisfying $M \models 2^{\kappa} = 2^{\kappa^+}$. Therefore, $U_\tau = \{ X \subset V_\kappa \mid \vec{U}\uhr \tau \in j(X)\}$ for all $\tau < \kappa^+$.
	Clearly, \[M \models 2^\kappa = 2^{\kappa^+} =  |([V_\kappa]^{<\omega} \times \power(V_\kappa))^{\kappa \times \kappa^+}|.\]
	In $V$, let $H : \kappa \to V_\kappa$ be a partial function with \[\dom(H) = \{\alpha < \kappa \mid 2^\alpha = 2^{\alpha^+}\},\] such that for every $\alpha \in \dom(H)$, $H(\alpha) : 2^\alpha \leftrightarrow ([V_\alpha]^{<\omega} \times \power(V_\alpha))^{\alpha \times \alpha^+}$ is a bijection.
	Let $H(\kappa) = j(H)(\kappa)$ (by slightly abuse of notations). By elementarity, $H(\kappa)$ is a bijection between $(2^\kappa)^M$ and $\left(([V_\kappa]^{<\omega} \times \power(V_\kappa))^{\kappa \times \kappa^+}\right)^M$. Note that both the domain and co-domain of $H(\kappa)$, are evaluated in $M$. It is likely to assume $H(\kappa)$ is not surjective on $\left(([V_\kappa]^{<\omega} \times \power(V_\kappa))^{\kappa \times \kappa^+}\right)^V$, as we do not assume that $M$ agrees with $V$ on $\mathcal{P}(\kappa^{+})$.
	
	For notational simplicity, we denote $H(\alpha)$ by $H_\alpha$, for each $\alpha \leq \kappa$.
	
	Let $R(\vec{U})$ be the Radin forcing associated with $\vec{U}$. Choose a generic set $G\subseteq R(\vec{U})$.
	We claim there is no weak diamond on $\kappa$ in $V[G]$.
	To show this, it will be convenient to identify  conditions  $p = \vec{d} \fr \la \vec{U},A\ra \in R(\vec{U})$ with pairs $\la \vec{d}, A\ra \in V_\kappa^{<\omega} \times \power(V_\kappa)$.
	We say that $\la \vec{d}, A\ra \in V_\kappa^{<\omega} \times \power(V_\kappa)$ is a simple representation of $p$.
	Since $R(\vec{U})$ satisfies $\kappa^+.c.c$ we can represent antichains in $R(\vec{U})$ using elements in $(V_\kappa^{<\omega} \times \power(V_\kappa))^\kappa$. We use the phrase simple representation in this context as well. Note that those antichains may have cardinality strictly less than $\kappa$ (and can even be empty). We code those antichains as well by appending the list of simple representations of the elements in the antichains with the empty set that does not represent any condition in $R(\vec{U})$.
	
	We first work in $V$. Let $\name{g}$ be an $R(\vec{U})$ name for a function from $\kappa$ to $2$ and $p = \vec{d} \fr \la \vec{U}, A\ra$ a condition.
	
	Without loss of generality (by taking a direct extension) we may assume that either $\len(\stem (p)) = 0$ or that for every $\vec\mu\in A$, $\kappa(\vec\mu) > \kappa_{\len(\stem (p)) - 1}$.
	
	Fix for a moment $\vec{\mu} \in A$ and consider the one-point extension \[p \fr \la \vec{\mu} \ra = \vec{d} \fr \la \vec{\mu}, A \cap V_{\kappa(\vec{\mu})}\ra \fr \la \vec{U}, A\ra.\] The forcing $R(\vec{U})/ (p \fr \la \vec{\mu}\ra)$ factors into a product \[R(\vec{\mu}) \times R(\vec{U})/ (\la \vec{U}, A \setminus V_{\kappa(\vec{\mu}) + 1}\ra).\]
	
	Let $\langle q_\xi \mid \xi < 2^{\kappa(\vec\mu)}\rangle$ be an enumeration of the conditions in $R(\vec{\mu})$. Let us define by induction sets $A^{\vec{\mu}}_\xi \in \bigcap \vec{U}$ such that $q_\xi \fr \langle \vec{U}, A^{\vec{\mu}}_\xi\rangle$ decides whether $\name{g}(\can{\kappa}(\vec{\mu})$ is $0$, $1$ or if $q_\xi$ does not force its value in the generic extension by the right component in the decomposition of $R(\vec{U}) / p$. Since $\cap \vec{U}$ is $\kappa$-complete and since $2^{\kappa(\vec{\mu})} < \kappa$, $A_{\vec\mu} = \bigcap_{\xi < 2^{\kappa(\vec{\mu}}} A^{\mu}_\xi \in \cap \vec{U}$.
	
	Clearly, for densely many conditions $q\in R(\mu)$,   $\la\vec{U}, A_{\vec\mu}\ra$ forces that they decide the value of $\name{g}(\can{\kappa}(\vec{\mu}))$.
	
	As $\vec\mu$ was arbitrary, by the considerations above, for every $\vec\mu \in A$, the condition $p \fr \vec{\mu}$ has a direct extension of the form
	$\vec{d} \fr \la \vec{\mu}, A \cap V_{\kappa(\vec{\mu})}\ra \fr \la \vec{U} , A_{\vec{\mu}}\ra$ forcing $\name{g}(\can{\kappa}(\vec{\mu})) = \sigma(\vec{\mu})$,
	where $\sigma(\vec{\mu})$ is a $R(\vec{\mu})$ name for an ordinal in $\{0,1\}$. Let \[A^* = \Delta_{\vec{\mu} \in A} A_{\vec{\mu}} = \{ \vec{\nu} \in V_\kappa \mid \vec{\nu} \in A_{\vec{\mu}} \text{ if } \vec{\mu} \in V_{\kappa(\vec{\nu})}\}\] and $p^* = \vec{d} \fr \la \vec{U}, A^*\ra$.  Thus, $p^*$ is a direct extension of $p$ and $p \fr \la \vec{\mu}\ra \force \name{g}(\can{\kappa}(\vec{\mu})) = \sigma(\vec{\mu})$ for all $\vec{\mu} \in A^*$.
	
	Let us fix a well order in $M$ of $V_{\kappa + \omega}$, $\trianglelefteq$.
	In $M$, we define a function $h : \kappa^+ \to  (V_\kappa^{<\omega} \times \power(V_\kappa))^\kappa$.
	
	For every $\tau < \kappa^+$, let $h(\tau) \in (V_\kappa^{<\omega} \times \power(V_\kappa))^\kappa$ be the $\trianglelefteq$-least simple representation of a maximal collection $A_\tau \subset R(\vec{U}\uhr \tau)$ of incompatible conditions
	$q \in R(\vec{U}\uhr \tau)$ which force $j(\sigma)(\vec{U}\uhr\tau) = \can{0}$.
	We point out the definition of $h$ relies on the assumption $\vec{U} \in M$, which implies $M$ contains a enumeration of the posets $\la R(\vec{U}\uhr \tau) \mid \tau < \kappa^+\ra$.
	With this enumeration in $M$, we can identify $h$ with an element in $\left(V_\kappa^{<\omega} \times \power(V_\kappa)^{\kappa \times \kappa^+}\right)^M$.
	
	Let $f = H_\kappa^{-1}(h) : \kappa \to 2$. $f$ is going to code $\name{g}$: we are going to define a name for a function $\name{F}\colon 2^{<\kappa} \to 2$ (in $V[G]$) such that $p^\star \Vdash \{\alpha < \check{\kappa}\mid \name{F}(f \restriction \alpha) = \name{g}(\alpha)\}$ contains a club.
	
	Back to $V$, let us define an auxiliary function $F' : 2^{<\kappa} \to V_\kappa$ as follows.
	
	Let $\alpha \in \dom H$, and $w\in 2^\alpha$. $F'(w)$ is a function whose domain is the collection of all measure sequences $\vec{\nu}$ with $\kappa(\vec{\nu}) = \alpha$. Let us look at $H_\alpha(w) \in (V_\alpha^{<\omega} \times \power(V_\alpha))^{\alpha \times \alpha^+}$. Let $\vec{\nu}$ be a measure sequence with $\kappa(\vec{\nu}) = \alpha$. We define $F'(w)(\vec{\nu})$ to be the set of all conditions $q \in R(\vec{\nu})$ such that $q$ is simply represented by an element of $H_\alpha(w)(\len(\vec{\nu}))$, where $\len(\vec{\nu})$ is the length of the sequence $\vec{\nu}$.
	
	Finally, we define $\name{F}$, which is forced to be a function from $2^{<\kappa}$ to $2$ in $V[G]$. For every $\vec{\nu} \in MS_G$, and $w \in 2^{\kappa(\vec{\nu})}$ we set
	$$\name{F}(w) =
	\begin{cases}
	0 &\mbox{ if } F'(w)(\vec{\nu}) \cap G(\vec{\nu}) \neq \emptyset\\
	1 &\mbox{ otherwise.}
	\end{cases}$$
	Define $\name{F}(w) = 0$ for every other $w \in 2^\alpha$ or for any $w\in 2^\alpha$ such that $\alpha \notin \dom(H)$. It is important to note that the definition of $\name{F}$ does not depend on $p$ or $f$.
	
	Let us show that $p^\star$ forces that the collection of all $\alpha < \kappa$ such that $\name{F}(f\restriction \alpha) = \name{g}(\alpha)$ contains a club. In fact, we will show that the condition $p^\star$ (which was defined above) forces that a tail of the set of accumulation point of the Radin club will be there.
	
	Let us work in $M$. $j(H)_\kappa(j(f) \uhr \kappa) = H_\kappa(f) = h$. Therefore, $j(F')(f)$ is a function, and its domain is the set of measure sequences with critical point $\kappa$. For $\gamma > 0$, $j(F')(f)(\vec{U} \uhr \gamma)$ is the collection of all conditions in $R(\vec{U} \uhr \gamma)$, which are simply represented by an element from $h(\gamma)$. By the definition of $h(\gamma)$ - this is a maximal set of incompatible conditions in $R(\vec{U} \uhr \gamma)$ which forces $j(\name{g})(\can{\kappa}) = \can{0}$.
	
	Let $X$ be the set of all $ \vec{\nu}\in MS_\kappa$ such that either $\len \vec{\nu} = 0$ or for $\alpha = \kappa(\nu)$, \[p^\star \fr \langle \vec{\nu}\rangle \Vdash F'(f\uhr \alpha)(\vec{\nu}) \cap G(\vec{\nu}) \neq \emptyset \iff \name{g}(\alpha) = \can{0}\]
	Then
	\[\vec{U} \uhr \gamma \in j(X)\]
	for all $\gamma < \kappa^{+}$. Therefore, $X\in \bigcap \vec{U}$.
	
	By elementarity, if $\vec\nu\in X$ and $\len \nu > 0$, then the condition $p^\star \fr \la \vec{\nu}\ra$, forces $\name{F}(f \uhr \kappa(\vec\nu)) = \name{g}(\kappa(\vec\nu))$.
	
	Let $\rho = \min \{\kappa(\vec{\nu} \mid \nu \in X\}$. We conclude that the condition $p^{\star\star} = \vec{d} \fr \la \vec{U}, A^\star \cap X\ra$ forces that for every $\alpha \in \text{acc } C_G \setminus \rho$, $\name{F}(f\restriction \alpha) = \name{g}(\alpha)$.
	\hfill \qedref{mt}
	\section{Weak diamond at weakly inaccessible cardinals}
	To round out the picture, we are left with the case of weakly but not strongly inaccessible cardinal. We distinguish three cases. If $2^{<\lambda}=2^\lambda$ then we already know that $\neg\Phi_\lambda$ and if $2^{<\lambda}=2^\kappa< 2^\lambda$ for some $\kappa<\lambda$ then $\Phi_\lambda$. The remaining case is when the sequence $\langle 2^\theta:\theta<\lambda\rangle$ is not eventually constant. In this case $2^{<\lambda}<2^\lambda$, and we do not know if the weak diamond holds (see \cite{MR2210145}, Question 1.28) though we have seen that the very weak diamond holds.
	It is possible to force $\Phi_\lambda$ in such cases:
	
	\begin{theorem}
		\label{posclm} It is consistent relative to the existence of an inaccessible cardinal that there is a weakly inaccessible, $\lambda$, such that $\langle 2^\theta:\theta<\lambda\rangle$ is not eventually constant, $\lambda$ is not strongly inaccessible and $\Phi_\lambda$ holds.
	\end{theorem}
	
	\par\noindent\emph{Proof}. \newline
	We begin with a model of ${\rm GCH}$ with a strongly inaccessible cardinal $\lambda$. We aim to blow up $2^\theta$ for every regular uncountable $\theta\leq\lambda$.
	
	Let $\rm Add(\mu, \nu)$ denote the standard forcing notion for adding $\nu$ many new subsets to $\mu$. Namely, it is the collection of all partial functions from $\mu \times \nu$ to $2$ with domain of size ${<}\mu$, ordered by inclusion. This forcing is $\mu$-complete and $(2^{{<}\mu})^{+}$-cc. (in particular, assuming $\rm GCH$, it is $\mu^{+}$-cc.). For basic properties of $\rm Add(\mu, \nu)$, see \cite[Chapter 15]{MR1940513}.
	
	Let $\mathbb{P}$ be the Easton support product of ${\rm Add}(\theta,\lambda^{+\theta+1})$ for every regular uncountable $\theta\leq\lambda$. Namely, it is the collection of all elements $p$ in the product
	\[\prod_{\theta \leq \lambda, \text{ regular}} {\rm Add}(\theta, \lambda^{+\theta + 1})\]
	such that ${\rm supp } (p) = \{\alpha \leq \lambda \mid p(\alpha) \neq \emptyset\}$ is an Easton set, i.e.\ $|{\rm supp} (p) \cap \rho| < \rho$ for all regular $\rho \leq \lambda$. Note that it $\lambda$ is the first inaccessible cardinal then Easton support is the same as $<\lambda$-support.
	
	Notice that $\mathbb{P}$ neither collapses cardinals nor changes cofinalities.
	
	The forcing $\mathbb{P}$ is $\aleph_1$-complete, as a product of $\aleph_1$-complete forcing notions with Easton support.
	
	For technical reasons, it will be more convenient to represent the conditions in $\mathbb{P}$ at partial functions from $\lambda^{+\lambda + 1}$ to $2$. Such a partial function is a condition if its support $s$ satisfies:
	\begin{enumerate}
		\item for every regular uncountable cardinal $\theta\leq \lambda$, $|s \cap [\lambda^{+\theta}, \lambda^{+\theta + 1})| < \theta$.
		\item the collection $t = \{\theta \leq \lambda \mid s \cap [\lambda^{\theta}, \lambda^{+\theta + 1}) \neq \emptyset\}$ consists only of regular cardinals and it is an Easton set, namely $|t \cap \alpha| < \alpha$ for every inaccessible cardinal $\alpha$.
	\end{enumerate}
	In particular, $|s| < \lambda$. From this point we will always assume that our conditions are represented in this way.
	
	Let us show that $\mathbb{P}$ is $\lambda{^+}$-cc. Let $\{p_i \mid i < \lambda^{+}\}$ be a collection of conditions in $\mathbb{P}$. Let $s_i \subseteq \lambda^{+|lambda + 1}$ be the support of $p_i$. Then by the $\Delta$-system lemma, there is a subcollection $I \subseteq \lambda^{+}$, and a set $r \subseteq \lambda^{+\lambda + 1}$, such that for every $i, j$ in $I$, if $i \neq j$ then $s_i \cap s_j = r$. Let us narrow down $I$ further to a set $J$ for which there is some constant function $q_\star$ such that for every $i$, $p_i \restriction r = q_\star$. This is possible since $2^{{<}\lambda} = \lambda < \lambda^{+}$ and $|r| < \lambda$. Thus, every pair of conditions $p_i, p_j$, such that $i,j\in J$ are compatible.
	
	By Easton's Theorem, $\langle 2^\theta:\theta<\lambda\rangle=\langle \lambda^{+\theta+1} : \theta<\lambda\rangle$ in the generic extension, so this sequence is not eventually constant. Although $\lambda$ is not strongly inaccessible any more, it is still weakly inaccessible. It remains to show that $\Phi_\lambda$ holds in $V^{\mathbb{P}}$.
	
	Let $\name{c}$ be a name of a coloring function from ${}^{<\lambda}2$ into $2$ in the generic extension. Observe that $2^{<\lambda} = \lambda^{+\lambda} <2^\lambda=\lambda^{+\lambda+1}$ in the generic extension $V^{\mathbb{P}}$.
	
	Let $\mathbb{P}'$ be the Easton support product of ${\rm Add}(\theta, \lambda^{+\theta + 1})$ over all $\theta < \lambda$. $\mathbb{P}$ is isomorphic to $\mathbb{P}' \times {\rm Add}(\lambda, \lambda^{+\lambda + 1})$. We claim that ${\rm Add}(\lambda, \lambda^{+\lambda + 1})^V$ is $\lambda$-distributive in the generic extension by $\mathbb{P}'$ and in particular, every element in $({}^{<\lambda}2 )^{V^{\mathbb{P}}}$ belongs to $V^{\mathbb{P}'}$. Let $t$ be a function from $\rho$ to the ordinals in the generic extension. Let us split the product into two components:
	
	\[\mathbb{P} \cong \Big(\prod_{\theta \leq \rho^{+}, \theta\text{ regular}}^{E} {\rm Add}(\theta, \lambda^{+\theta + 1})\Big) \times \Big(\prod_{\rho^{++} \leq \theta \leq \lambda,\,\theta\text{ regular}}^{E} {\rm Add}(\theta, \lambda^{+\theta + 1})\Big)\]
	where the $E$ denotes that the product is with Easton support. The first component is $\rho^{++}$-cc (by the same argument as the chain condition of $\mathbb{P}$) and the second argument is $\rho^{++}$-complete. Thus, by Easton's lemma, the second component is $\rho^{++}$-distributive in the generic extension by the first component. In particular, it cannot add $t$.
	
	Let us pick, in $V$, for every potential member of ${}^{<\lambda}2$ in the generic extension, a name $\name{t} \subseteq \mathbb{P}' \times \lambda$. Every member of ${}^{<\lambda}2$ will be equal to an evaluation of one of those names by the distributivity of ${\rm Add}(\lambda, \lambda^{+\lambda + 1}$. There are at most $(\lambda^{+\lambda})^\lambda = \lambda^{+\lambda}$ such names. For each name $\name{t}$ as above, let us pick a maximal antichain of conditions in $\mathbb{P}$ that decides the value of $\name{c}(\name{t})$. Let $B$ be the union of the supports of all the conditions that appear in one of those antichains. Recall that the support of a condition in $\mathbb{P}$ is a subset of $\lambda^{+\lambda + 1}$. Thus, $B$ is the union of $\lambda^{+\lambda}$ sets of size $<\lambda$ and thus $|B| = \lambda^{+\lambda}$, and in particular, it is bounded. Let $\delta = \sup B < \lambda^{+\lambda + 1}$.
	
	For a set $C$, let us denote by $\mathbb{P}\upharpoonright C$ the set of conditions $p\in \mathbb{P}$ such that $\mathrm{supp}\ p \subseteq C$. Clearly, $\mathbb{P}\cong \mathbb{P}\upharpoonright B\times \mathbb{P}\upharpoonright (\lambda^{+\lambda+1}\setminus B)$ and $\name{c}\in V^{\mathbb{P}\upharpoonright B}$. Denote $\mathbb{P}\upharpoonright B$ by $\mathbb{Q}$ and $\mathbb{P}\upharpoonright (\lambda^{+\lambda+1}\setminus B)$ by $\mathbb{R}$. Since $\mathbb{Q}$ is $\aleph_1$-complete, in the generic extension by $\mathbb{Q}$, $\mathbb{R}$ is still $\aleph_1$-complete.
	
	Let $\name{G}$ be the canonical name for the generic filter of $\mathbb{P}$. $\bigcup \name{G}$ is forced to be a function from $\lambda^{+\lambda + 1}$ to $2$. Let $\name{g}$ be the canonical name for the $\delta$-th new $\lambda$-Cohen set, namely for $\alpha < \lambda$, $p \Vdash \alpha \in \name{g}$ iff $p(\lambda^{+\lambda} + \lambda \cdot \delta + \alpha) = 1$.
	
	Let us work in $V^{\mathbb{Q}}$,
	Let $\name{f}$ be an $\mathbb{R}$-name of a function from $\lambda$ into $2$, and let $\name{D}$ be an $\mathbb{R}$-name of a club subset of $\lambda$. For every condition $r\in\mathbb{R}$ we can choose by recursion a sequence of conditions $\langle r_n:n\in\omega\rangle$ in $\mathbb{R}$ and a sequence of ordinals $\la \beta_n \mid n < \omega\ra$:
	\begin{enumerate}
		\item [$(a)$] $r_0=r$ and $r_n\leq r_{n+1}$.
		\item [$(b)$] $r_{n+1}\Vdash \name{f}\upharpoonright \alpha_n = \check{g}_n$ for some $g_n\in V^{\mathbb{Q}}$.
		\item [$(c)$] $\{\xi \mid \lambda^{+\lambda} + \lambda \cdot \delta + \xi \in \dom r_{n+1}\} = \alpha_n$
		\item [$(d)$] $r_{n+1}\Vdash \check{\beta}_n\in \name{D}$ and $\check{\alpha}_n<\check{\beta}_n<\check{\alpha}_{n+1}$.
	\end{enumerate}
	Let $p$ be $\bigcup\limits_{n\in\omega}r_n$. By the $\aleph_1$-completeness of $\mathbb{R}$, $p$ is a condition in $\mathbb{R}$. Denote $\bigcup\limits_{n\in\omega}\alpha_n$ by $\alpha$. Since $r_{n+1}\leq p$ for every $n\in\omega$ we see that $p\Vdash\name{f}\upharpoonright\alpha= \bigcup\limits_{n\in\omega}g_n$. Likewise, $p\Vdash\check{\alpha}\in\name{D}$ since $\name{D}$ is forced to be closed and $\alpha=\bigcup\limits_{n\in\omega}\beta_n$.
	
	Our goal is to show that the fixed $g$ chosen above can serve as a weak diamond function from $\lambda$ into $2$. For this, we shall prove that the condition $p$ can be extended to force $g(\alpha)= \name{c}(\name{f}\upharpoonright \check{\alpha})$.
	As $\name{c}\in V^{\mathbb{Q}}$, the value of $\name{c}(\name{f}\upharpoonright \check{\alpha})$ is determined by the condition $p$. However, functions from $\lambda$ into $2$ are not determined in a bounded stage. In particular, we can extend $p$ to a condition $q$ which forces $g(\alpha)= \name{c}(\name{f}\upharpoonright \check{\alpha})$.
	It follows that for every $\name{f}:\check{\lambda}\rightarrow 2$ and every club $\name{D}$ there exists an ordinal $\check{\alpha}\in\name{D}$ for which $\Vdash_{\mathbb{P}} g(\alpha)=c(f\upharpoonright\alpha)$, so we are done.
	
	\hfill \qedref{posclm}
	\section{Conclusions and Open Problems}
	The following diagram summarizes the relationship between the various prediction principles considered in this paper:
	
	
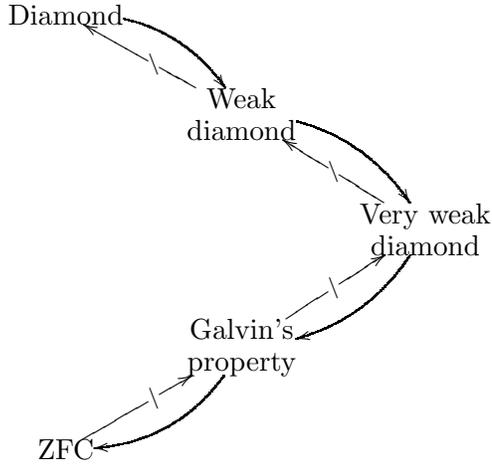
\xymatrix{
		*\txt{Diamond} \ar @/^1pc/ [dr] \\
		& *\txt{Weak\\ diamond} \ar @/^1pc/ [dr] \ar @{->}[ul]|\setminus \\
		& & *\txt{Very weak\\ diamond} \ar @/^1pc/ [dl] \ar @{->}[ul]|\setminus \\
		& *\txt{Galvin's\\ property} \ar @/^1pc/ [dl] \ar @{->}[ur]|\setminus \\
		*\txt{ZFC} \ar @{->}[ur]|\setminus
	}
	
	\medskip
	
	The downward positive implications $\Diamond_\lambda\Rightarrow\Phi_\lambda\Rightarrow\Psi_\lambda$ are trivial. The fact that $\Psi_\lambda$ implies Galvin's property appears in \cite{MR3604115}, as well as the negative direction upwards (i.e., Galvin's property does not imply $\Psi_\lambda$). The consistency of $\Psi_\lambda$ with $\neg\Phi_\lambda$ can be exemplified by a strongly inaccessible cardinal for which $\neg\Phi_\lambda$ is forced. The consistency of $\Phi_\lambda$ with $\neg\Diamond_\lambda$ can be forced by simple cardinal arithmetic considerations. Finally, it is shown in \cite{MR830084} that Galvin's property is not a theorem of ZFC, as its negation can be forced.
	
	We mention the question of Shelah from \cite{MR2210145}, which seems to be the last open case:
	
	\begin{question}
		\label{q1} Assume $\lambda$ is weakly inaccessible and $\langle 2^\theta:\theta<\lambda\rangle$ is not eventually constant. Is it consistent that $\neg\Phi_\lambda$ holds?
	\end{question}
	
	\newpage
	\providecommand{\bysame}{\leavevmode\hbox to3em{\hrulefill}\thinspace}
	\providecommand{\MR}{\relax\ifhmode\unskip\space\fi MR }
	\providecommand{\MRhref}[2]{%
		\href{http://www.ams.org/mathscinet-getitem?mr=#1}{#2}
	}
	\providecommand{\href}[2]{#2}
	
\end{document}